\documentclass[12pt]{article}
\usepackage{amsmath,amssymb,amsthm}
\newtheorem{theorem}{Theorem}

\renewcommand{\le}{\leqslant}
\renewcommand{\ge}{\geqslant}

\DeclareMathOperator{\rank}{rank}

\begin{document}

\title{Canonical form of
$m$-by-$2$-by-$2$ matrices over a
field of characteristic other than two\footnotetext{This is a preliminary version of the paper published in Linear Algebra Appl. 418 (2006) 15--19.}}
\author{Genrich Belitskii%
\thanks{Partially supported by
Israel Science Foundation, Grant
186/01.}\\ Dept. of Mathematics,
Ben-Gurion University of the Negev\\
Beer-Sheva 84105, Israel,
 genrich@cs.bgu.ac.il
 \and
M. Bershadsky,\\ Sapir Academic College P.b.\\ Hof
Ashkelon, 79165, Israel,
maximb@mail.sapir.ac.il
       \and
Vladimir V. Sergeichuk%
\thanks{
The research was done while this author
was visiting the Ben-Gurion University
of the Negev.}\\
Institute of Mathematics,
Tereshchenkivska 3, Kiev,
Ukraine\\sergeich@imath.kiev.ua}
\date{}

\maketitle

\begin{abstract}
We give a canonical form of $m\times 2\times 2$ matrices for equivalence over any field of characteristic not two.

{\it AMS classification:} 15A21; 15A69

{\it Keywords:} Multidimensional matrices;
Tensors; Classification
 \end{abstract}

Each trilinear form $f: U\times
V\times W\to \mathbb F$ on vector
spaces $U$, $V$, and $W$ of
dimensions $m$, $n$, and $q$ over a
field $\mathbb F$ is given by the
$m\times n\times q$ spatial matrix
\begin{equation}\label{1.0}
{\cal A} =[a_{ijk}]_{i=1}
^m{}_{j=1}^n{}_{k=1}^q,\qquad
a_{ijk}=f(u_i,v_j,w_k),
\end{equation}
with respect to bases $\{u_i\}$ of
$U$, $\{v_j\}$ of $V$, and $\{w_k\}$
of $W$. Changing the bases, we can
reduce ${\cal A}$ by {\it
equivalence transformations}
\begin{equation}\label{1.1}
[a_{ijk}]\longmapsto
[b_{i'j'k'}],\qquad
b_{i'j'k'}=\sum_{ijk}
a_{ijk}r_{ii'}s_{jj'} t_{kk'},
\end{equation}
in which $[r_{ii'}]$, $[s_{jj'}]$,
and $[t_{kk'}]$ are nonsingular
matrices.

Canonical forms of $2\times 2\times
2$ matrices for equivalence over the
field of complex numbers was given
by Ehrenborg \cite{2}. We extend his
result to $m\times 2\times 2$
matrices over any field $\mathbb F$
of characteristic not two. Unlike
\cite{2}, one of our canonical
$2\times 2\times 2$ matrices depends
on a parameter that is determined up
to multiplication by any $z^2,\ 0\ne
z\in\mathbb F$.

Note that the canonical form problem
for $m\times n\times 3$ matrices for
equivalence is hopeless since it
contains the problem of classifying
pairs of linear operators and, hence,
the problem of classifying arbitrary
systems of linear operators
\cite[Theorems 4.5 and 2.1]{1}; such
classification problems are called
\emph{wild}. The canonical form problem
for $m\times n\times 2$ matrices for
equivalence was studied in
\cite[Section 4.1]{1}.

We give the spatial matrix
\eqref{1.0} by the $q$-tuple of
$m\times n$ matrices
\begin{equation}\label{1.3}
{\cal A}=\|A_1\,|\dots|\, A_q\|,
\qquad A_k=[a_{ijk}]_{ij}.
\end{equation}
The equivalence transformation
\eqref{1.1} can be realized in two
steps: by the nonsingular linear
substitution
\begin{equation}\label{1.5}
C_1=A_1t_{11}+\dots+A_qt_{q1},\
\dots,\
C_q=A_1t_{1q}+\dots+A_qt_{qq},
\end{equation}
and then by the simultaneous
equivalence transformation
\begin{equation}\label{1.4}
\|R^TC_1S\,|\dots|\,R^TC_qS\|,\qquad
R=[r_{ii'}],\ S=[s_{jj'}].
\end{equation}
Moreover, two spatial matrices are
equivalent if and only if one
reduces to the other by nonsingular
linear substitutions and
simultaneous equivalence
transformations of the form
\eqref{1.5} and \eqref{1.4}. The
rank
\[
q'=\rank \{A_1,\dots,A_q\}
\]
of the matrices $A_1,\dots,A_q$ from
\eqref{1.3} in the space of
$m$-by-$n$ matrices is an invariant
of ${\cal A}$ with respect to
equivalence transformations.

Apart from \eqref{1.3}, the spatial
matrix \eqref{1.0} can be also given
by the tuples
\begin{equation}\label{a20}
(\tilde{A_1},\dots,
\tilde{A_n}),\qquad (\Tilde{\Tilde
A}_1,\dots, \Tilde{\Tilde A}_m)
\end{equation}
of the matrices
$\tilde{A_j}=[a_{ijk}]_{ik}$ and
$\Tilde{\Tilde A}_i=[a_{ijk}]_{jk}$,
and
\begin{equation}\label{a3}
n'=\rank \{\tilde A_1,\dots,\tilde
A_n\},\qquad m'=\rank \{\Tilde{\Tilde
A}_1,\dots,\Tilde{\Tilde A}_m\}
\end{equation}
are also invariants of ${\cal A}$
for equivalence transformations. We
say that the spatial matrix ${\cal
A}$ is \emph{regular} if $m=m'$,
$n=n'$, and $q=q'$.

Let ${\cal A}$ be irregular. Make
the first $q'$ matrices
$A_1,\dots,A_{q'}$ in the $q$-tuple
\eqref{1.3} linearly independent and
the others $A_{q'+1},\dots,A_{q}$
zero by substitutions of the form
\eqref{1.5}. Then reduce in the same
way the tuples \eqref{a20} of the
obtained spatial matrix and get a
spatial matrix ${\cal B} =[b_{ijk}]$
that is equivalent to $\cal A$ and
whose entries outside of
\begin{equation}\label{a6}
{\cal B}'
=[b_{ijk}]_{i=1}^{m'}{}_{j=1}^{n'}
{}_{k=1}^{q'}
\end{equation}
are zero. We call $\cal B$ a
\emph{regularized form of} $\cal A$
and its submatrix \eqref{a6} a {\it
regular part of} $\cal A$.

Two spatial matrices of the same size are equivalent if and only if their regular parts are equivalent
\cite[Lemma 4.7]{1}. Hence, it suffices to give a canonical form of a regular spatial matrix.

Due to this lemma and the next
theorem, the set of all regularized
forms whose regular submatrices are
\eqref{a9}--\eqref{a14z} is a set of
$m\times 2\times 2$ canonical
matrices for equivalence.

\begin{theorem}\label{t.1}
Over a field\/ $\mathbb F$ of
characteristic not two, each regular
$m\times n\times q$ matrix $\cal A$
with $n\le 2$ and $q\le 2$ is
equivalent to one of the spatial
matrices:
\begin{equation}\label{a9}
\begin{Vmatrix}
  \,1\,
\end{Vmatrix}
\qquad(1\times 1\times 1),
\end{equation}
\begin{equation}\label{a10}
\begin{Vmatrix}
 \, 1&0\,\\0&1
\end{Vmatrix}
\qquad(2\times 2\times 1),
\end{equation}
\begin{equation}\label{a10a}
\left|\!\left|\!\begin{array}{c|c}
    1&0\\ 0&1
\end{array}\!\right|\!\right|
\qquad(2\times 1\times 2),
\end{equation}
\begin{equation}\label{a10b}
\left|\!\left|\!\begin{array}{cc|cc}
    1&0&0&1
\end{array}\!\right|\!\right|
\qquad(1\times 2\times 2),
\end{equation}
\begin{equation}\label{a11}
\|I_2\,|\,B(a)\|:=\left|\!\left|\!\begin{array}{cc|cc}
    1&0&0& a\\ 0&1&1&0
\end{array}\!\right|\!\right|
\qquad(a\in\mathbb F,\
2\times 2\times 2),
\end{equation}
\begin{equation}\label{a12}
\left|\!\left|\!\begin{array}{cc|cc}
    1&0&0&0 \\ 0&1&0&0\\ 0&0&0&1
\end{array}\!\right|\!\right|
\qquad(3\times 2\times 2),
\end{equation}
\begin{equation}\label{a12a}
\left|\!\left|\!\begin{array}{cc|cc}
    1&0&0&0 \\ 0&1&1&0\\ 0&0&0&1
\end{array}\!\right|\!\right|
\qquad(3\times 2\times 2),
\end{equation}
\begin{equation}\label{a14z}
\left|\!\left|\!\begin{array}{cc|cc}
    1&0&0&0 \\ 0&1&0&0\\
    0&0&1&0\\0&0&0&1
\end{array}\!\right|\!\right|
\qquad(4\times 2\times 2).
\end{equation}
These spatial matrices, except for
\eqref{a11}, are uniquely determined by
$\cal A$; $\|I_2\,|\,B(a)\|$ is
equivalent to $\|I_2\,|\,B(b)\|$ if and
only if $a=bz^2$ for some nonzero
$z\in\mathbb F$.
\end{theorem}

\begin{proof} \emph{Existence.}
Let ${\cal A}$ be $m\times n\times
1$. Since ${\cal A}=\|A\|$ is
regular, it reduces by elementary
transformations \eqref{1.4} to
\eqref{a9} or \eqref{a10}.

Let ${\cal A}$ be $1\times 2\times
2$. Reduce it to the form
$\|1\,0\,|\,b_1\,b_2\|$ by
transformations \eqref{1.4}. Since
$\cal A$ is regular, $b_2\ne 0$; we
make $b_2=1$ multiplying the second
columns by $b_1^{-1}$. Adding a
multiple of $b_2$, make $b_1=0$ and
obtain \eqref{a10b}. Similarly, if
${\cal A}$ is $2\times 1\times 2$,
then it reduces to \eqref{a10a}.
Note that \eqref{a10a} and
\eqref{a10b} can be obtained from
the spatial matrix $[a_{ijk}]$
defined in \eqref{a10} by
interchanging its indices.

It remains to consider ${\cal A}$ of
size $m\times 2\times 2$ with $m\ge
2$. Since ${\cal A}=\|A\,|\,B\|$ is
regular, $A\ne 0$, $B\ne 0$, and the
rows of the $m\times 4$ matrix $[A\;
B]$ are linearly independent; that
is, $m =\rank [A\; B]\le 4.$
Interchanging $A$ and $B$ if
necessary, we make
\begin{equation}\label{a30}
\rank A\ge\rank B.
\end{equation}

If $m =4$, then $[A\; B]$ is a
nonsingular $4\times 4$ matrix and
we reduce ${\cal A}=\|A\,|\,B\|$ by
row-transformations to the form
\eqref{a14z}.

If $m= \rank [A\; B]=3$, then $\rank
A=2$ by \eqref{a30}, we reduce
${\cal A}=\|A\,|\,B\|$ by
transformations \eqref{1.4} to the
form
\begin{equation}\label{a17}
\left|\!\left|\!\begin{array}{cc|cc}
    1&0&b_{11}&0 \\
    0&1&b_{21}&0\\ 0&0&0&1
\end{array}\!\right|\!\right|.
\end{equation}
Replacing $B$ by $B-b_{11}A$, we
make $b_{11}=0$ but spoil
$b_{22}=0$; we fix $b_{22}$ adding
the third row. So $\cal A$ is
equivalent to \eqref{a12} or
\eqref{a12a}.

Let $m= \rank [A\; B]=2$. If $A$ is
singular, then by \eqref{a30} $\rank
A=\rank B=1$ and we reduce ${\cal
A}$ to the form
\begin{equation*}\label{a15zz}
\left|\!\left|\!\begin{array}{cc|cc}
  1 & 0 &0&b_{12}\\
    0 & 0&b_{21}&b_{22}
\end{array}\!\right|\!\right|.
\end{equation*}
Since $\cal A$ is regular,
$(b_{21},b_{22})$ and
$(b_{12},b_{22})$ are nonzero.
Because $\rank B=1$, $b_{21}= 0$ or
$b_{12}= 0$. We replace $A$ by the
nonsingular $A+B$.

Since $A$ is nonsingular, $\cal A$
reduces to the form
\begin{equation}\label{a19}
\left|\!\left|\!\begin{array}{cc|cc}
  1 & 0 &b_{11}&b_{12}\\
    0 & 1&b_{21}&b_{22}
\end{array}\!\right|\!\right|.
\end{equation}
Preserving $A=I_2$, we will reduce
$B$ by similarity transformations.
If $b_{21}=0$, then we make
$b_{21}\ne 0$ using
\[
\begin{bmatrix}
  1& 0\\ -\varepsilon &1
\end{bmatrix}
\begin{bmatrix}
b_{11}&b_{12}\\ b_{21}&b_{22}
\end{bmatrix}
\begin{bmatrix}
  1& 0\\ \varepsilon &1
\end{bmatrix},\qquad\text{$\varepsilon=1$
or $\varepsilon=-1$.}
\]
Multiplying the first row of $B$ by
$b_{21}$ and its first column by
$b_{21}^{-1}$, we obtain $b_{21}=1$.
Make $b_{11}=-b_{22}$ replacing $B$
by $B-\alpha A= B-\alpha I_2$, where
$\alpha := (b_{11}+b_{22})/2$. At
last, replace $B$ with
\[\begin{bmatrix}
  1& -b_{11}\\ 0&1
\end{bmatrix}
\begin{bmatrix}
b_{11}&b_{12}\\ 1&-b_{11}
\end{bmatrix}
\begin{bmatrix}
  1& b_{11}\\ 0&1
\end{bmatrix}
=\begin{bmatrix} 0&\star\\ 1&0
\end{bmatrix}
\]
and obtain \eqref{a11}.
\medskip

\noindent \emph{Uniqueness.} If two
spatial matrices among
\eqref{a9}--\eqref{a14z} are
equivalent, then they have the same
size, and so they are \eqref{a12}
and \eqref{a12a}, or
$\|I_2\,|\,B(a)\|$ and
$\|I_2\,|\,B(b)\|$. If the spatial
matrix \eqref{a12} is $\|A\,|\,B\|$,
then $\rank({\alpha A+\beta B})\ne
1$ for all $\alpha,\beta\in\mathbb
F$, hence \eqref{a12} is not
equivalent to \eqref{a12a}.

Let $\|I_2\,|\,B(a)\|$ be equivalent
to $\|I_2\,|\,B(b)\|$. By
\eqref{1.4} and \eqref{1.5}, there
is a nonsingular matrix $
\left[\begin{smallmatrix} \alpha
&\beta
\\ \gamma &\delta
\end{smallmatrix}\right]
$ such that the matrices $\alpha
I_2+\beta B(a)$ and $\gamma
I_2+\delta B(a)$ are simultaneously
equivalent to $I_2$ and $B(b)$. Then
$\alpha I_2+\beta B(a)$ is
nonsingular and the matrix $(\alpha
I_2+\beta B(a))^{-1}(\gamma
I_2+\delta B(a))$ is similar to
$B(b)$. Hence, the matrices
\begin{equation*}\label{a20x}
\begin{bmatrix}
\alpha &-\beta a \\ -\beta &\alpha
\end{bmatrix}
\begin{bmatrix}
\gamma &\delta a\\ \delta &\gamma
\end{bmatrix}=
\begin{bmatrix}
\alpha\gamma-\beta\delta a
&(\alpha\delta -\beta \gamma)a
\\ \alpha\delta -\beta \gamma
&\alpha\gamma-\beta\delta a
\end{bmatrix},\quad
(\alpha^2-\beta^2
a)^2\begin{bmatrix} 0 &b \\ 1 &0
\end{bmatrix}
\end{equation*}
are similar. Equating their traces
and determinants we obtain
$\alpha\delta -\beta \gamma=0$ and $
(\alpha\delta -\beta \gamma)^2a=
(\alpha^2-\beta^2 a)^2 b. $
Therefore, $a=bz^2$, where
$z=(\alpha\delta -\beta \gamma)/
(\alpha^2-\beta^2 a)$.

Conversely, if $a=bz^2$ and $0\ne
z\in\mathbb F$, then
$\|I_2\,|\,B(b)\|$ is equivalent to
$\|I_2\,|\,zB(b)\|$, which is
equivalent to $\|I_2\,|\,B(a)\|$
since
\[
\begin{bmatrix}
  z& 0\\ 0&1
\end{bmatrix}
\begin{bmatrix}
0&bz\\ z&0
\end{bmatrix}
\begin{bmatrix}
  z^{-1}& 0\\ 0&1
\end{bmatrix}
=\begin{bmatrix} 0&bz^2\\ 1&0
\end{bmatrix}.
\]
\end{proof}

\end{document}